\documentclass[12pt, notitlepage]{article}
\makeindex
\usepackage{amsfonts,amssymb,amsmath,amscd,mathtools,amsthm,xcolor,mathrsfs,tikz-cd}
\usepackage[backref,colorlinks, linkcolor=blue, citecolor=black]{hyperref}
\usepackage{latexsym}
\usepackage[top=1in, bottom=1.25in, left=1.25in, right=1.25in]{geometry}
\makeatletter

\newtheorem{thm}{Theorem}[section]
\newtheorem{prop}[thm]{Proposition}

\newtheorem{lem}[thm]{Lemma}

\numberwithin{equation}{section}

\def\Nrd{{\rm Nrd}}

\def\Gal{{\rm Gal}}
\def\Ind{{\rm ind\,}}

\title{On triviality of the reduced Whitehead group over Henselian fields}
\author{Abhay Soman}
\date{}
\usepackage{subfiles}
\begin{document}
 \maketitle
\begin{abstract}
	Let $F$ be a Henselian field of $q$-cohomological dimension $3$, where $q$ is a prime. Let $\Gamma_F$ be the totally ordered abelian value group of $F$ and let $D$ be a central division algebra over $F$ of index a power of $q$ such that the characteristic of the residue field, $\overline{F}$ is coprime to $q$. We show that when $1\leq\dim_{\mathbb{F}_q}(\Gamma_F/q\Gamma_F)\leq 3$, the reduced Whitehead group of $D$ is trivial.
\end{abstract}

\section{Introduction}
Let $E$ be an arbitrary field and $A$ be a finite-dimensional central simple algebra over $E$. We denote the group of units of $A$ by $A^*$. The \emph{reduced Whitehead group} of $A$ is given by \[SK_1(A)=\{a\in A^*: \Nrd_A(a)=1\}/[A^*,A^*],\] where $[A^*,A^*]$ is the commutator subgroup of $A^*$ and $\Nrd_A$ is the reduced norm map, $\Nrd_A:A^*\rightarrow E^*$. If $A\simeq M_n(D)$ then, there is an isomorphism $SK_1(A)\simeq SK_1(D)$ (cf. \S $23$, corollary $1$, \cite{D}). Moreover, if $D_i$ ($1\leq i\leq r$) are central division algebras of $p_i$-power degrees, where $p_i$ are distinct primes then, $SK_1(D_1\otimes D_2\otimes\cdots\otimes D_r)\simeq SK_1(D_1)\times SK_1(D_2)\times\cdots\times SK_1(D_r)$ (cf. \S $23$, lemma $6$, \cite{D}). Thus, to study $SK_1$ it is enough to consider central division algebras of prime power degrees.

Let $G_E$ denote the absolute Galois group of a field $E$, i.e., $G_E=\Gal(E^s/E)$, where $E^s$ is the separable closure of $E$. The $q$-\emph{cohomological dimension} of $E$ is the least positive integer $d$ such that for all discrete $G_E$-modules $A$ which are $q$-primary torsion groups, the Galois cohomology groups $H^i(G_E,A)$ are trivial for $i\geq d+1$. The $q$-cohomological dimension of $E$ is denoted by $cd_q(E)$.
\medskip

For a torsion free abelian group $\Gamma_F$ and a prime $q$, let the $q$-rank of $\Gamma_F$ be $r_q\coloneqq\dim_{\mathbb{F}q}(\Gamma_F/q\Gamma_F)$. In this article our aim is to prove the following theorem.

\begin{thm}\label{intro main}
	Let $(F,v)$ be a Henselian field with totally ordered abelian value group $\Gamma_F$ and characteristic of the residue field $\overline{F}$, $char(\overline{F})=\bar{p}$. Let $q\neq\bar{p}$ be a prime and $D$ be a finite-dimensional central division algebra over $F$ of index a power of $q$. If  $cd_q(F)=3$ and $1\leq r_q\leq 3$ then, $SK_1(D)=(1)$.
\end{thm}

In \cite{S2}, Suslin conjectured that $SK_1(D)=(1)$ for any central division algebra $D$ of index $q^2$ where $q$ is a prime, over fields of cohomological dimension $3$. The above theorem provides further evidence for the validity of the conjecture.

In order to prove the theorem \ref{intro main} we first obtain a relation between the $q$-cohomological dimension of a Henselian field and its residue field. More precisely, we prove the following proposition. 
\begin{prop}\label{intro_cd}
	Let $(F,v)$ be a Henselian field with totally ordered abelian value group $\Gamma_F$ and characteristic of the residue field $\overline{F}$, $char(\overline{F})=\bar{p}$. For a prime $q\neq \bar{p}$, assume that $r_q$ is finite. Then we have
	\[cd_q(F)=cd_q(\overline{F})+r_q.\]
\end{prop}

Along with the above proposition (\ref{intro_cd}), the proof of theorem (\ref{intro main}) uses valuation theory on division algebras over Henselian fields as developed in \cite{Ershov} and \cite{TW-book}.
\section{Preliminaries}
\subsection{Notations and basic definitions}
In this section we recall some notions on division algebras over Henselian fields. We refer the reader to \cite{TW-book} for more details.

Let $(F,v)$ be a Henselian valued field with the value group $v(F^*)=\Gamma_F$. We assume that $\Gamma_F$ is a totally ordered additive abelian group. All division algebras in this article are assumed to be finite-dimensional over its center. Let $D$ be a central division algebra over $F$. The valuation $v$ extends uniquely to the valuation $v_D$ on $D$ (cf. corollary $1.7$, \cite{TW-book}) and by (theorem $1.4$, \cite{TW-book}) it is given by
\[
v_D(d):=\frac{1}{\Ind D} v(\Nrd_D(d)), \text{ for $d\in D^*$}.
\] We denote the value group of $v_D$ by $\Gamma_D=v_D(D^*)$. We denote the \emph{residue division algebra} by $\overline{D}=\{x\in D: v_D(x)\geq 0 \}/\{x\in D: v_D(x)>0 \}$, and its center by $Z(\overline{D})$. There is a well-defined group homomorphism (cf. \S $1.1.1$, \cite{TW-book})
\[
\theta_D:\Gamma_D\rightarrow\textrm{Aut}(Z(\overline{D})/\overline{F}),
\]
where $\overline{F}$ is the residue field of $F$ with respect to the valuation $v$. We call $\theta_D$ the \emph{canonical homomorphism} of the valuation $v_D$.

If $D$ satisfies following equality \[[D:F]=[\overline{D}:\overline{F}]\,|\Gamma_D:\Gamma_F|,\] then $D$ is said to be \emph{defectless} over $F$. The division algebra $D$ is said to be \emph{tame} over $F$ if $D$ is defectless over $F$, the extension $Z(\overline{D})/\overline{F}$ is separable and $char(\overline{F})\nmid |\ker(\theta_D):\Gamma_F|$. Recall that $D$ is said to be \emph{totally ramified} over $F$ if $|\Gamma_D:\Gamma_F|=[D:F]$ or equivalently, $\overline{D}=\overline{F}$ and $D$ is defectless (cf. \S $7.4.1$, \cite{TW-book}).

For a \emph{tame} central division algebra $D$, we denote by $G$ the \emph{abelian Galois} group $\Gal(Z(\overline{D})/\overline{F})$ (cf. proposition $1.5$, \cite{TW-book}).
Let \[\widetilde{N}=N_{Z(\overline{D})/\overline{F}}\circ\Nrd_{\overline{D}}:\overline{D}\rightarrow \overline{F}, \quad\text{and}\quad \zeta=\Ind D/(\Ind\overline{D})\, [Z(\overline{D}):\overline{F}].
\]
If we further assume that $\overline{D}$ is a field, then $\zeta^2=[D:F]/[\overline{D}:\overline{F}]^2$. By (proposition $1.5$, \cite{TW-book}), $\overline{D}/\overline{F}$ is an abelian Galois extension, and there is an isomorphism of abelian groups $\Gal(\overline{D}/\overline{F})\simeq\Gamma_D/\ker\theta_D$. Furthermore, as $D$ is defectless 
\begin{equation}\label{value of zeta}
\zeta^2=|\Gamma_D:\ker\theta_D||\ker\theta_D:\Gamma_F|/[\overline{D}:\overline{F}]=|\ker\theta_D:\Gamma_F|.
\end{equation}
We associate following groups with the valuation $v_D$ (cf. Ershov \cite{Ershov}).
\begin{description}
	\item $SL(D)=\{x\in D^*:\Nrd_D(x)=1\}$;
	\item $U=\{x\in D^*:v_D(x)=0 \}$;
	\item $SL^{v_D}(D)=\{x\in SL(D):\widetilde{N}(\overline{x})=1 \}$;
	\item $SK^{v_D}(D)=SL^{v_D}(D)/[U,D^*]$;
	\item $\mathcal{K}=\left(\overline{[U,D^*]}\cap SL(\overline{D})\right)\big/[\overline{D}^*,\overline{D}^*]$;
	\item $C_{v_D}=\mu_{\zeta}(\overline{F})\cap\widetilde{N}(\overline{D}^*)$.
\end{description}

We recall the definition of the $-1$-Tate cohomology group. Let $G$ be a finite abelian group and $M$ be a $G$-module, with operation written multiplicatively. Let $N_G\colon M\to M$ be the $G$-norm map given by $m\mapsto\prod_{\sigma\in G}\sigma(m)$. Let $I_G(M)=\langle m\sigma(m)^{-1}:m\in M,\sigma\in G\rangle$ be a subgroup of $M$. Then, the $-1$-Tate cohomology group of $G$ with respect to $M$ is \[\widehat{H}^{-1}(G,M)=\ker N_G/I_G(M).\]

\subsection{Some known results}
We recall some known results which are used in the proof of main theorem. We refer to the earlier subsection for notations and terminology. We have the following important result due to Ershov.
\begin{thm}[page $68$, \cite{Ershov}]
	Let $F$ be a Henselian field and let $D$ be a tame central division algebra over $F$. With notations as in the above subsection, we have the following diagram with exact rows and column.
	
	\begin{center}\label{diagram}
		\begin{tikzcd}
		& 1 \arrow[d]\\
		1 \to \mathcal{K} \to SK_1(\overline{D}) \arrow[r] & SK^{v_D}(D) \arrow[r] \arrow[d] & \widehat{H}^{-1}(G,\Nrd_{\overline{D}}(\overline{D}^*)) \to 1\\	
		
		1\to [D^*,D^*]/[U,D^*] \arrow[r] & SL(D)/[U,D^*] \arrow[d] \arrow[r] & SK_1(D) \to 1\\
		& C_{v_D}\arrow[d]\\
		& 1
		\end{tikzcd}
	\end{center}
	
\end{thm}

Now we state a result on the triviality of the reduced Whitehead group of division algebras of $q$-prime power indices over a field $K$ of $q$-cohomological dimension at most $2$. More precisely,
\begin{thm}(theorem $1.1$, \cite{Boley})\label{Boley}
	Let $K$ be a field of characteristic $p$ (which can be zero) and let $q$ be a prime number different from $p$. Suppose that $cd_q(K)\leq 2$ and that $A$ is a central simple algebra over $K$ whose index is a power of $q$. Then, the reduced Whitehead group of $A$ is trivial.
\end{thm}

\section{Computation of cohomological dimension}\label{Computation of cohomological dimension}
Let $F$ be a Henselian field with valuation $v$ and $char(\overline{F})=\overline{p}$. We denote by $F^t$ the inertia field of $F$ and $G^t=\Gal(F^s/F^t)$, its corresponding Galois group (see \S $5.2$, \cite{EP}). We denote by $G_F$ the absolute Galois group of the field $F$. By (theorem $5.2.7$, \cite{EP}) the absolute Galois group of the residue field $\overline{F}$, $G_{\overline{F}}\simeq\Gal(F^t/F)\simeq G_{F}/G^t$. Consider the inflation map induced by the above isomorphism: 
\[
\textrm{inf}: H^i(G_{\overline{F}},\mathbb{Z}/n)\rightarrow H^i(G_{F},\mathbb{Z}/n); \text{ for } (n,\overline{p})=1.
\]
By (theorem $5.3.3$ ($1$), \cite{EP}), the ramification Galois group $G^v=\Gal(F^s/F^v)$ is the unique $\overline{p}$-Sylow subgroup of $G^t$. For $\bar{p}\nmid n$, we thus have, $H^i(G^v,\mathbb{Z}/n)=0$ for $i\geq 1$ (cf. I.\S$3.3$, corollary $2$, \cite{Se1}). Thus, by the inflation-restriction sequence (VII.\S$6$, proposition $5$, \cite{Se}) we get the following isomorphism 

\[
H^i(G^t/G^v, \mathbb{Z}/n)\simeq H^i(G^t,\mathbb{Z}/n), \text{ for } i\geq 0 \text{ and } (n,\bar{p})=1.
\] 
\medskip

Let $q$ be a prime distinct from $\overline{p}$. By (theorem $5.3.3$ ($3$), \cite{EP}), there is an isomorphism of profinite groups 
\begin{equation}\label{G^t/G^v}
G^t/G^v\simeq\prod_{q\neq\overline{p}} \mathbb{Z}_q^{r_q}
\end{equation}
where for each prime $q\neq\overline{p}$, $r_q$ is the $q$-rank of $\Gamma_F$, i.e., $r_q$ is the $\mathbb{F}_q$-dimension of $\Gamma_F/q\Gamma_F$. Thus, for $i\geq 0$ and gcd$(n,\bar{p})=1$ we have

\begin{equation}\label{coho of G^t}
H^i(G^t,\mathbb{Z}/n)\simeq H^i(G^t/G^v,\mathbb{Z}/n)\simeq H^i(\prod_{q\neq\overline{p}}\mathbb{Z}_q^{r_q},\mathbb{Z}/n).
\end{equation}
\medskip

We first prove the following lemma which appears as a part of an exercise in \cite{Se1}.
\begin{lem}[I.\S $4.5$, Exercise $(1)$, \cite{Se1}]\label{Z_q}
	For a prime number $q$ and a natural number $r_q$ we have, \[cd_q(\mathbb{Z}_q^{r_q})= r_q \text{ and } H^{r_q}(\mathbb{Z}_q^{r_q},\mathbb{Z}/q) \text{ is a group of order $q$.}\]
\end{lem}
\begin{proof}
	We prove this lemma by induction on $r_q$. Suppose $r_q=1$, then $cd_q(\mathbb{Z}_q)=1$ and $H^1(\mathbb{Z}_q,\mathbb{Z}/q)$ is a group of order $q$. We assume the result for $r_q-1<\infty$, i.e., $cd_q(\mathbb{Z}_q^{r_q-1})=r_q-1$ and $H^{r_q-1}(\mathbb{Z}_q^{r_q-1},\mathbb{Z}/q)$ is a group of order $q$. Now to prove the result for $r_q$, consider the short exact sequence
	\[1\rightarrow \mathbb{Z}_q\rightarrow \mathbb{Z}_q^{r_q}\rightarrow \mathbb{Z}_q^{r_q-1}\rightarrow 1.\]
	By the induction hypothesis $cd_q(\mathbb{Z}_q^{r_q-1})=cd_q(\mathbb{Z}_q^{r_q}/\mathbb{Z}_q)=r_q-1<\infty$ and $cd_q(\mathbb{Z}_q)=1$. Moreover, $\mathbb{Z}_q$ is a pro-$q$-group and $H^1(\mathbb{Z}_q,\mathbb{Z}/q)$ is a group of order $q$. Therefore, by (I.\S $4.1$, proposition $22$, \cite{Se1})
	\[
	cd_q(\mathbb{Z}_q^{r_q})=r_q.
	\]
	By the spectral sequence (I.\S $3.3$, remark, \cite{Se1}),
	\[
	H^{r_q}(\mathbb{Z}_q^{r_q},\mathbb{Z}/q)=H^{r_q-1}(\mathbb{Z}_q^{r_q-1}, H^1(\mathbb{Z}_q,\mathbb{Z}/q))=H^{r_q-1}(\mathbb{Z}_q^{r_q-1},\mathbb{Z}/q).
	\] which is a group of order $q$ by the induction hypothesis. Hence the lemma is proved.
\end{proof}
\medskip

We now prove the proposition \ref{intro_cd}. The proof is on similar lines as the discrete valued case (cf. II.\S$4.3$, proposition $12$, \cite{Se1}) and we present it here for completeness.
\medskip

\noindent{\it Proof of proposition} \ref{intro_cd}:

Consider the following exact sequence 
\[
1\rightarrow H\rightarrow \prod_{q'\neq \overline{p}}\mathbb{Z}_{q'}^{r_{q'}}\rightarrow\mathbb{Z}_q^{r_q}\rightarrow 1.
\]
where the product is taken over all primes $q'\neq\overline{p}$ and the last map is projection onto the $q$th component. By (I.\S $3.3$, corollary $2$, \cite{Se1}), $cd_q(H)=0$ and by (lemma \ref{Z_q}) $cd_q(\mathbb{Z}_q^{r_q})=r_q < \infty$.  Since the group $\prod_{q'\neq\overline{p}}\mathbb{Z}_{q'}^{r_{q'}}$ is abelian, using (I.\S $4.1$, proposition $22$, \cite{Se1}), we get
\begin{equation}\label{cd_r_q}
cd_q(\prod_{q'\neq\overline{p}}\mathbb{Z}_{q'}^{r_{q'}})=r_q.
\end{equation}
In fact, using inflation-restriction sequence (VII.\S$6$, proposition $5$, \cite{Se}) we get an isomorphism
\begin{equation}\label{G^t as a order q group}
H^{r_q}(\mathbb{Z}_q^{r_q},\mathbb{Z}/q)\simeq H^{r_q}(\prod_{q'\neq\overline{p}}\mathbb{Z}_{q'}^{r_{q'}},\mathbb{Z}/q).
\end{equation} In particular, $H^{r_q}(\prod_{q'\neq\overline{p}}\mathbb{Z}_{q'}^{r_{q'}},\mathbb{Z}/q)$ is a group of order $q$ (see lemma \ref{Z_q}). Hence, $H^{r_q}(G^t,\mathbb{Z}/q)$ is also a group of order $q$ (see equation \ref{coho of G^t}), and applying (I.\S 3.3, proposition $15$, \cite{Se1}) to the groups $G^t, G^v$ we also have $cd_q(G^t)\leq r_q$. As a result, $cd_q(G^t)=r_q$. Moreover, using equation (\ref{cd_r_q}) above and (I.\S $3.3$, proposition $15$, \cite{Se1}) we have  
\[
cd_q(G_F)\leq cd_q(G_{\overline{F}})+r_q.
\]
Assume $cd_q(G_F/G^t)(=cd_q(G_{\overline{F}}))=n<\infty$, using the spectral sequence (I.\S $3.3$, remark, \cite{Se1}), we get 
\[
H^{n+r_q}(G_F,\mathbb{Z}/q)=H^n(G_F/G^t, H^{r_q}(G^t,\mathbb{Z}/q)).
\]

Since by our assumption $cd_q(G_F/G^t)=n$ and as observed above, $H^{r_q}(G^t,\mathbb{Z}/q)$ is a group of order $q$, we get $H^n(G_F/G^t,H^{r_q}(G^t,\mathbb{Z}/q))\neq 0$. Hence, $cd_q(G_F)=cd_q(G_{\overline{F}})+r_q$.

\section{Triviality of the reduced Whitehead group}
%\subsection{Lemmata}
In this section we give a proof of theorem \ref{intro main}. We start with the following lemma. We denote the absolute Galois group of a field $E$ by $G_E$.
\begin{lem}\label{cohom trivial module}
	Let $L/E$ be a  Galois field extension of degree $q^r$ for some $r\in\mathbb{N}$ and a prime number $q$. Further assume that the $q$-cohomological dimension of $E$, $cd_q(G_E)\leq 1$. Then the multiplicative $Gal(L/E)$-module $L^*$ is cohomologically trivial.
\end{lem}
\begin{proof}
	We show that $H^1(\Gal(L/E),L^*)$ and $H^2(\Gal(L/E),L^*)$ are trivial. Indeed, by Hilbert's theorem $90$, $H^1(\Gal(L/E),L^*)=0$ (cf. X.\S 1, proposition $2$, \cite{Se}). By (X.\S 4, corollary, \cite{Se}) the $q$-group, $H^2(\Gal(L/E), L^*)$ injects into the Brauer group of $E$. Since $cd_q(G_E)\leq 1$, the $q$-primary torsion part of the Brauer group is trivial. Hence $H^2(\Gal(L/E),L^*)=0$. Thus, the Tate {\it cohomology groups with positive exponents} $1$ and $2$ are trivial.  Using (IX.\S 5, theorem $8$, \cite{Se}) we get, the $\Gal(L/E)$-module $L^*$ is cohomologically trivial, i.e., for every subgroup $H\leq \Gal(L/E)$ and every integer $i$, $\widehat{H}^i(H,L^*)=0$.
\end{proof}

A. R. Wadsworth suggested the proof of the following statement.
\begin{lem}\label{cyclicity}
	Let $\Gamma_F\subset\Gamma_D$ be totally ordered abelian groups with $\Gamma_F$ as a subgroup of $\Gamma_D$. Suppose that the cardinality of $\Gamma_D/\Gamma_F$, $|\Gamma_D/\Gamma_F|=q^n$, where $q$ is a prime number, and $n$ is a natural number. Assume that $r_q$ is finite. Then we have 
	\begin{enumerate}
		\item If $r_q=0$ then, $\Gamma_D=\Gamma_F$.
		\item If $r_q=1$ then, $\Gamma_D/\Gamma_F$ is cyclic. 
	\end{enumerate}
	In general, the number of invariant factors of $\Gamma_D/\Gamma_F$ is at most $r_q$.
\end{lem}	
\begin{proof}
	We have $\Gamma_D/\Gamma_F\simeq q^n\Gamma_D/q^n\Gamma_F\subset \Gamma_F/q^n\Gamma_F$. Furthermore, as $\Gamma_F$ is a torsion-free abelian group, the multiplication by $q$ map from  $q^{i-1}\Gamma_F/q^{i}\Gamma_F$ to $q^{i}\Gamma_F/q^{i+1}\Gamma_F$ is a group isomorphism for every $i\geq 1$. Hence we have
	\[
	\Gamma_F/q\Gamma_F\simeq q\Gamma_F/q^2\Gamma_F\simeq\ldots\simeq q^i\Gamma_F/q^{i+1}\Gamma_F\simeq\ldots\simeq q^{n-1}\Gamma_F/q^n\Gamma_F.
	\]
	Thus, if $r_q=0$, then $\Gamma_F=q\Gamma_F=\ldots=q^n\Gamma_F$; hence $\Gamma_D=\Gamma_F$. If $r_q=1$, the above isomorphism shows $|\Gamma_F/q^n\Gamma_F|=q^n$. The finite abelian $q$-group $\Gamma_F/q^n\Gamma_F$ is a direct product of cyclic groups of $q$-power order; it is therefore cyclic, since $(\Gamma_F/q^n\Gamma_F)/q(\Gamma_F/q^n\Gamma_F)\simeq\Gamma_F/q\Gamma_F$, which is cyclic. Therefore, $\Gamma_D/\Gamma_F$ is cyclic.
	
	As $\dim_{\mathbb{F}_q}(\Gamma_F/q\Gamma_F)=r_q$, the isomorphism $(\Gamma_F/q^n\Gamma_F)/q(\Gamma_F/q^n\Gamma_F)\simeq\Gamma_F/q\Gamma_F$ show that the number of invariant factors of $\Gamma_D/\Gamma_F$, which is a subgroup of $\Gamma_F/q^n\Gamma_F$, is at most $r_q$.
\end{proof}

We recall the definition of a symbol algebra appearing in the next lemma. Let $E$ be a field and $n\geq 2$ an integer such that $E$ contains a primitive $n$-th root of unity $\omega$. For any $a,b\in E^*$, consider the $E$-algebra $(a,b)_n$ generated by two elements $i,j$ subject to the relations $i^n=a,j^n=b, ij=\omega ji$. This algebra is called a symbol algebra.
\begin{lem}\label{totally ramified is symbol}
	We keep the notations of theorem \ref{intro main}. If $r_q=2$ or $3$ then, every nonsplit tame totally ramified central division algebra $T$ over $F$ of index a power of $q$ is a symbol algebra. Moreover, $SK_1(T)=(1)$.
\end{lem}
\begin{proof}
	Suppose that the degree of $T$ is $q^t$. By (proposition $7.72$ and corollary $7.76$, \cite{TW-book}), $\overline{F}$ contains a primitive root of unity of order $\exp(\Gamma_T/\Gamma_F)$ and the invariant factors of $\Gamma_T/\Gamma_F$ occur in pairs, and $\exp(T)=\exp(\Gamma_T/\Gamma_F)$. Since $T$ is a nonsplit tame totally ramified division algebra, $\Gamma_T\neq\Gamma_F$. Hence, $\Gamma_T/\Gamma_F\simeq\mathbb{Z}/q^t\mathbb{Z}\times\mathbb{Z}/q^t\mathbb{Z}$ by (lemma \ref{cyclicity}). Therefore, $\exp(T)=q^t=\Ind (T)$. By (theorem $11.23$ (ii), \cite{TW-book}), $SK_1(T)\simeq\mu_{\Ind(T)}(\overline{F})/\mu_{\exp(T)}(\overline{F})=(1)$. Furthermore, $T$ is a symbol algebra by (proposition $7.74$, \cite{TW-book}).
	
\end{proof}

We now proceed to prove the main theorem \ref{intro main} stated in the introduction.
\medskip

\noindent{\it Proof of theorem} \ref{intro main}:
\medskip

%The $q$-cohomological dimension of $\overline{F}$, $cd_q(\overline{F})<\infty$, since the absolute Galois group of $\overline{F}$, $G_{\overline{F}}$ is isomorphic to a subgroup of $\Gal(F^v/F)$, where $F^v$ is the ramification field for the extension $F_s/F$. 
By (proposition \ref{intro_cd}), $cd_q(F)=cd_q(\overline{F})+r_q$. Note that the division algebra $D$ over the Henselian field $F$ is tame. Indeed, we have \[[D:F]=\partial_{D/F}\cdot [\overline{D}:\overline{F}]\cdot |\Gamma_D:\Gamma_F|, \] where $\partial_{D/F}$ is the defect and $\partial_{D/F}=\bar{p}^l$ for some $l\in\mathbb{N}$ (cf. theorem $4.3$, \cite{TW-book}). Since, $q\neq\bar{p}$ and $D$ has the degree a power of $q$, $\bar{p}\nmid |\ker({\theta_D}):\Gamma_F|$ and the defect, $\partial_{D/F}=1$, i.e., in particular, $D$ is defectless. Similarly, $\bar{p}\nmid [Z(\overline{D}):\overline{F}]$ and hence $Z(\overline{D})/\overline{F}$ is a separable field extension. Moreover by (proposition $1.5$ (iii), \cite{TW-book}), the field extension $Z(\overline{D})/\overline{F}$ is normal. Hence, $Z(\overline{D})/\overline{F}$ is a Galois extension. Put $G=\Gal(Z(\overline{D})/\overline{F})$. 
\smallskip

We first consider the case $r_q=1$. In this case, $cd_q(\overline{F})=2$. As $Z(\overline{D})/\overline{F}$ is a finite extension of fields, we have $cd_q(Z(\overline{D}))= 2$ (see II.\S $4$, proposition $10$, \cite{Se1}). Thus, by (theorem \ref{Boley}) $SK_1(\overline{D})=(1)$.  The canonical homomorphism $\theta_D:\Gamma_D\rightarrow \Gal(Z(\overline{D})/\overline{F})$ is surjective by (proposition $1.5$, \cite{TW-book}). Clearly we have $\Gamma_F\subseteq\ker(\theta_D)\subseteq\Gamma_D$. By (lemma \ref{cyclicity}), $\Gamma_D/\Gamma_F$ is a cyclic group and hence $\Gamma_D/\ker(\theta_D)\simeq \Gal(Z(\overline{D})/\overline{F})$ is a cyclic group as well. By Hilbert's theorem $90$ (X.\S $1$, corollary, \cite{Se}), $\widehat{H}^{-1}(G,Z(\overline{D})^*)=(1)$. As $\deg(\overline{D})$ is a power of $q$ and $cd_q(Z(\overline{D}))\leq 2$, by Merkurjev-Suslin theorem (theorem $24.8$, \cite{S1}) $\Nrd_{\overline{D}}(\overline{D}^*)=Z(\overline{D})^*$. Hence, \[\widehat{H}^{-1}(G,\Nrd_{\overline{D}}(\overline{D}^*))=\widehat{H}^{-1}(G,Z(\overline{D})^*)=(1).\] Since $\ker\theta_D/\Gamma_F$ is a cyclic group, it can not carry a nondegenerate alternating pairing if it is not trivial. In view of (proposition $8.17$ (iv), \cite{TW-book}), we have $\ker\theta_D=\Gamma_F$. Hence, $\zeta=1$ and $SK_1(D)=(1)$ by (theorem \ref{diagram}).
\medskip

Now assume that $r_q=2$. We require the following result (cf. corollary (c), page $69$, \cite{Ershov})
\begin{equation}\label{exponent relation}
\exp(\Gamma_D/\Gamma_F)=\exp([D^*,D^*]/[U,D^*]).
\end{equation}
By (proposition $8.59$, \cite{TW-book}), $D\sim S\otimes_F T\in Br(F)$, where $S$ is an inertially split central division algebra and $T$ is a tame totally ramified central division algebra over $F$. Since $cd_q(\overline{F})\leq 1$, $Z(\overline{D})=\overline{D}$ and thus, $SK_1(\overline{D})=(1)$. By (lemma \ref{cohom trivial module}), $\widehat{H}^0(G,\overline{D}^*)\coloneqq\overline{F}^*/N_{\overline{D}/\overline{F}}(\overline{D}^*)=(1)$ and  $\widehat{H}^{-1}(G,\Nrd_{\overline{D}}(\overline{D}^*))=\widehat{H}^{-1}(G,\overline{D}^*)=(1)$. Therefore, $\widetilde{N}(\overline{D}^*)=\overline{F}^*$, and  $\mu_\zeta(\overline{F})\cap\widetilde{N}(\overline{D}^*)=\mu_\zeta(\overline{F})$.  Furthermore, exactness of first row in the diagram of (theorem \ref{diagram}) implies that
\begin{equation}\label{r_q=2 triviality of column}
SK^{v_D}(D)=(1) \quad\text{and}\quad SL(D)/[U, D^*]\simeq\mu_\zeta(\overline{F}).
\end{equation}

If $T= F$, then $D\sim S\in Br(F)$. Hence, $SK_1(D)\simeq SK_1(S) $. Thus, we may assume that $D$ is an inertially split central division algebra. In this case, $\zeta=1$ (cf. proposition $8.64$, \cite{TW-book}). So by (theorem \ref{diagram}) and (eq. \ref{r_q=2 triviality of column}), $SK_1(D)=(1)$. Now assume that $T\neq F$. %We first compute $\zeta$. Since $Z(\overline{D})=\overline{D}$, \[\zeta \coloneqq\Ind D/(\Ind\overline{D})\, [Z(\overline{D}):\overline{F}]=\Ind D/ [\overline{D}:\overline{F}].\] As $D$ is defectless, i.e., $[D:F]=[\overline{D}:\overline{F}]\, |\Gamma_D:\Gamma_F|$, \[\zeta^2=|\Gamma_D:\ker\theta_D||\ker\theta_D:\Gamma_F|/[\overline{D}:\overline{F}]. \] 
%There is a group isomorphism, $\Gamma_D/\ker\theta_D\simeq \Gal(\overline{D}/\overline{F})$ by (proposition $1.5$ (iv), \cite{TW-book}).
By (eq. \ref{value of zeta}), $\zeta^2=|\ker\theta_D:\Gamma_F|$. By (theorem $8.60$, \cite{TW-book}), $\exp(\ker\theta_D/\Gamma_F)=\exp(\Gamma_T/\Gamma_F)$, and $|\ker\theta_D:\Gamma_F|=|\Gamma_T:\Gamma_F|=[T:F]=q^{2t}$, where $\deg T=q^t$. Therefore, by (lemma \ref{totally ramified is symbol})
\begin{equation}\label{exponent inequality}
\zeta=q^t=\exp(\ker\theta_D/\Gamma_F)\leq\exp(\Gamma_D/\Gamma_F),
\end{equation}
and by (proposition $8.17$(v), \cite{TW-book}) $\zeta$-th primitive root of unity belongs to $\overline{F}$. The isomorphism  $SL(D)/[U, D^*]\simeq\mu_{\zeta}(\overline{F})$  implies that $SL(D)/[U, D^*]$ is a cyclic group of exponent $q^t$. Therefore, $\exp([D^*,D^*]/[U,D^*])\leq q^t$. By (eq. \ref{exponent relation} and eq. \ref{exponent inequality}), 
\begin{equation}\label{exp equals exp of totally ramified}
\exp(\Gamma_D/\Gamma_F)=\exp([D^*,D^*]/[U,D^*])=q^t.
\end{equation}
Therefore, $SL(D)/[U, D^*]=[D^*,D^*]/[U,D^*]$. In view of the following exact sequence (cf. theorem \ref{diagram}) \[1\rightarrow [D^*,D^*]/[U,D^*]\rightarrow SL(D)/[U, D^*]\rightarrow SK_1(D)\rightarrow 1, \]  we get that, $SK_1(D)=(1)$.
\medskip

If $r_q=3$ then, $cd_q(\overline{F})=0$ and the $q$-Sylow subgroup of the absolute Galois group of $\overline{F}$ is trivial (cf. I.\S$3.3$, corollary $2$, \cite{Se1}). Thus, $\overline{F}$ does not have a proper separable field extension of degree a power of $q$. For a defectless central division algebra $D$ over $F$ whose index is a power of $q$, we thus have $\overline{D}=Z(\overline{D})=\overline{F}$, and $[D:F]=|\Gamma_D:\Gamma_F|$, i.e., $D$ is totally ramified. Hence, by (lemma \ref{totally ramified is symbol}), $SK_1(D)=(1)$.

% ------------------------------------------------------------------------

\subsection*{Acknowledgment}
The author wishes to sincerely thank Jean-Pierre Tignol for carefully reading the manuscript and suggesting improvements. The author is grateful to A. R. Wadsworth for pointing out an error in the case $r_q=2$ of the theorem \ref{intro main}, in an earlier version of this manuscript. A. R. Wadsworth and J.-P. Tignol's critical comments immensely helped in the improvement of the manuscript. Thanks are also due to B. Sury, Amit Kulshrestha and Sudhanshu Shekhar for helpful discussions. We thank the reviewer for helpful suggestions.

\end{document}